\documentclass{article}
\usepackage{amssymb}
\usepackage{amsmath}

\newtheorem{theorem}{Theorem}

\newtheorem{example}[theorem]{Example}

\newtheorem{lemma}[theorem]{Lemma}

\newtheorem{proposition}[theorem]{Proposition}
\newtheorem{remark}[theorem]{Remark}

\input{tcilatex}
\begin{document}

\title{Minimum divergence estimators, Maximum Likelihood and the generalized
bootstrap}
\author{Michel Broniatowski \\
%EndAName
LPSM, CNRS\ UMR\ 8001, Sorbonne-Universit\'{e} Paris, France}
\maketitle

\begin{abstract}
This paper is an attempt to set a justification for making use of some

dicrepancy indexes, starting from the classical Maximum Likelihood
definition, and adapting the corresponding basic principle of inference to
situations where minimization of those indexes between a model and some
extension of the empirical measure of the data appears as its natural
extension. This leads to the so called generalized bootstrap setting for
which minimum divergence inference seems to replace Maximum Likelihood one.

\end{abstract}

%%%%%%%%%%%%%%%%%%%%%%%%%%%%%%%%%%%%%%%%%%

%%%%%%%%%%%%%%%%%%%%%%%%%%%%%%%%%%%%%%%%%%

\section{Motivation and context}

Divergences between probability measures are widely used in Statistics and
Data Science in order to perform inference \ under models of various kinds,
parametric or semi parametric, or even in non parametric settings.\ The
corresponding methods extend the likelihood paradigm and insert inference in
some minimum "distance" framing, which provides a convenient description for
the properties of the resulting estimators and tests, under the model or
under misspecification.\ Furthermore they pave the way to a large number of
competitive methods, which allows for trade-off between efficiency and
robustness, among others.\ Many families of such divergences have been
proposed, some of them stemming from classical statistics (such as the
Chi-square), while others have their origin in other fields such as
Information theory.\ Some measures of discrepancy involve regularity of the
corresponding probability measures while others seem to be restricted to
measures on finite or countable spaces, at least when using them as
inferential tools, henceforth in situations when the elements of a model
have to be confronted with a dataset. The choice of a specific discrepancy
measure in specific context is somehow arbitrary in many cases, although the
resulting conclusion of the inference might differ accordingly, above all
under misspecification; however the need for such approaches is clear when
aiming at robustness.

This paper considers a specific class of divergences, which contains most of
the classical inferential tools, and which is indexed by a single scalar
parameter.\ This class of divergences belongs to the
Csiszar-Ali-Silvey-Arimoto family of divergences (see \cite{LV1987}), and is
usually referred to as the power divergence class, which has been considered
By Cressie and Read \cite{CR1988}; however this denomination is also shared
by other discrepancy measures of some different nature\ \cite{Basu98}; see 
\cite{BrVaj} for a comprehensive description of those various inferential
tools with a discussion on their relations.\ We will use the acronym CR for
the class of divergences under consideration in this paper.

We have tried to set a justification for those discrepancy indexes, starting
from the classical Maximum Likelihood definition, and adapting the
corresponding basic principle of inference to situations where those indexes
appear as its natural extension. This leads to the so called generalized
bootstrap setting for which minimum divergence inference seems to replace
Maximum Likelihood one.

The contents of this approach can be summarized as follows.

Section 2 states that the MLE is obtained as a proxy of the minimizer of the
Kullback-Leibler divergence between the generic law of the observed variable
and the model, which is the large deviation limit for the empirical
distribution. This limit statement is nothing but the continuation of the
classical ML paradigm, namely to make the dataset more "probable " under the
fitted distribution in the model, or, equivalently, to fit the most "likely"
distribution in the model to the dataset.

Section 3 states that given a divergence pseudo distance $\phi $ in CR the
Minimum Divergence Estimator (MDE) is \ obtained as a proxy of the minimizer
of the large deviation limit for some bootstrap version of the empirical
distribution, which establishes that the MDE is MLE for bootstrapped samples
defined in relation with the divergence. This fact is based on the strong
relation which associates to any CR $\phi $ -divergence a specific RV\ $W$ \
(see Section \ref{Subsection Weights}) ; this link is the cornerstone for
the interpretation of the minimum $\phi $-divergence estimators as MLE's for
specific bootstrapped sampling schemes where $W$ has a prominent r\^{o}le.\
Some specific remark explores the link between MDE and MLE in exponential
families. As a by product we also introduce a bootstrapped estimator of the
divergence pseudo-distance $\phi $ between the distribution of the data and
the model.\ 

In Section 4 we specify the bootstrapped estimator of the divergence which
can be used in order to perform an optimal test of fit.\ Due to the type of
asymptotics handle in this paper, optimality is studied in terms of Bahadur
efficiency.\ It is shown that tests of fit based on such estimators enjoy
Bahadur optimality with respect to other bootstrap plans when the bootstrap
is performed under the distribution associated with the divergence criterion
itself.\ 

The discussion held in this paper pertains to parametric estimation in a
model $\mathcal{P}_{\Theta}$ whose elements $P_{\theta}$ are probability
measures defined on the same finite space $\mathcal{Y}:=\left\{
d_{1},..,d_{K}\right\} $, and $\theta\in\Theta$ an index space; we assume
identifiability, namely different values of $\theta$ induce different
probability laws $P_{\theta}$'s. Also all the entries of $P_{\theta}$ will
be positive for all $\theta$ in $\Theta.$

\subsection{Notation}

\subsubsection{Divergences}

We consider regular \textit{divergence functions} $\varphi$ which are non
negative convex functions with values in $\overline{\mathbb{R}^{+}}$ which
belong to $C^{2}$ $\left( \mathbb{R}\right) $ and satisfy $\varphi\left(
1\right) =\varphi^{\prime}\left( 1\right) =0$ and $\varphi^{\prime\prime
}\left( 1\right) =1$\textit{; }see \cite{LV1987} and \cite{BrSt2019} for
properties and extensions \textit{.} An important class of such functions is
defined through the power divergence functions 
\begin{equation}
\varphi_{\gamma}\left( x\right) :=\frac{x^{\gamma}-\gamma x+\gamma-1}{%
\gamma\left( \gamma-1\right) }  \label{CRdiv}
\end{equation}
defined for all real $\gamma\neq0,1$ with $\varphi_{0}\left( x\right)
:=-\log x+x-1$ (the likelihood divergence function) and $\varphi_{1}\left(
x\right) :=x\log x-x+1$ (the Kullback-Leibler divergence function). This
class is usually referred to as the Cressie-Read family of divergence
functions (see \cite{CR1988}). It is a very simple class of functions (with
the limits in $\gamma\rightarrow0,1$) which allows to represent nearly all
commonly used statistical criterions.\ Parametric inference in commonly met
situations including continuous models or some non regular models can be
performed with them; see\cite{BK2009}. The $L_{1}$ divergence function $%
\varphi\left( x\right) :=\left\vert x-1\right\vert $ is not captured by the
CR family of functions. When undefined the function $\varphi$ is declared to
assume value +$\infty.$

Associated with a divergence function $\varphi ,$ $\phi $ is the \textit{%
divergence pseudo-distance} between a probability measure and a finite
signed measure; see \cite{BrVaj}.

For $P:=\left( p_{1},..,p_{K}\right) $ and $Q:=\left( q_{1},..,q_{K}\right) $
in $\mathbb{S}^{K},$ the simplex of all probability measures on $\mathcal{Y}$%
, define, whenever $Q$ and $P$ have non null entries%
\begin{equation*}
\phi \left( Q,P\right) :=\sum_{k=1}^{K}p_{k}\varphi \left( \frac{q_{k}}{p_{k}%
}\right) .
\end{equation*}%
Indexing this pseudo-distance by $\gamma $ and using $\varphi _{\gamma }$ as
divergence function yields the likelihood divergence $\phi _{0}\left(
Q,P\right) :=-\sum p_{k}\log \left( \frac{q_{k}}{p_{k}}\right) $, the
Kullback-Leibler divergence $\phi _{1}\left( Q,P\right) :=\sum q_{k}\log
\left( \frac{q_{k}}{p_{k}}\right) $, the Hellinger divergence $\phi
_{1/2}\left( Q,P\right) :=\frac{1}{2}\sum p_{k}\left( \sqrt{\left( \frac{%
q_{k}}{p_{k}}\right) }-1\right) ^{2}$, the modified (or Neyman) $\chi ^{2}$
divergence $\phi _{-1}\left( Q,P\right) :=\frac{1}{2}\sum p_{k}\left( \left( 
\frac{q_{k}}{p_{k}}\right) -1\right) ^{2}\left( \frac{q_{k}}{p_{k}}\right)
^{-1}$. The $\ \chi ^{2}$ divergence $\phi _{2}\left( Q,P\right) :=\frac{1}{2%
}\sum p_{k}\left( \left( \frac{q_{k}}{p_{k}}\right) -1\right) ^{2}$ is
defined between signed measures; see \cite{BK2006} for definitions in more
general setting, and \cite{BK2009} for the advantage to extend the
definition to possibly signed measures in the context of parametric
inference for non regular models. Also the present discussion which is
restricted to finite spaces $\mathcal{Y}$ can be extended to general spaces.

The conjugate divergence function of $\varphi$ is defined through%
\begin{equation}
\widetilde{\varphi}\left( x\right) :=x\varphi\left( \frac{1}{x}\right)
\label{conjugatediv}
\end{equation}
and the corresponding divergence pseudo-distance $\widetilde{\phi}\left(
P,Q\right) $ is 
\begin{equation*}
\widetilde{\phi}\left( P,Q\right) :=\sum_{k=1}^{K}q_{k}\widetilde{\varphi }%
\left( \frac{p_{k}}{q_{k}}\right)
\end{equation*}
which satisfies%
\begin{equation*}
\widetilde{\phi}\left( P,Q\right) =\phi\left( Q,P\right)
\end{equation*}
whenever defined, and equals $+\infty$ otherwise. When $\varphi=\varphi
_{\gamma}$ then $\widetilde{\varphi}=\varphi_{1-\gamma}$ as follows by
substitution. Pairs $\left( \varphi_{\gamma},\varphi_{1-\gamma}\right) $ are
therefore \textit{conjugate pairs}. Inside the Cressie-Read family, the
Hellinger divergence function is self-conjugate.

For $P=P_{\theta}$ and $Q\in\mathbb{S}^{K}$ we denote $\phi\left( Q,P\right) 
$ by $\phi\left( Q,\theta\right) $ (resp $\phi\left( \theta,Q\right) $, or $%
\phi\left( \theta^{\prime},\theta\right) ,$ etc according to the context).

\subsubsection{Weights\label{Subsection Weights}}

This paragraph introduces the special link which connects CR divergences
with specific random variables, which we call weights. Those will be
associated to the dataset and define what is usually referred to as a
generalized bootstrap procedure.\ This is the setting which allows for an
interpretation of the MDE's as  generalized bootstrapped  MLE's.

For a given real valued random variable (RV) $W$ denote%
\begin{equation}
M(t):=\log E\exp tW  \label{mgfW}
\end{equation}%
its cumulant generating function which we assume to be finite in a non void
neighborhood of $0$ $.$ The Fenchel Legendre transform of $M$ \ (also called
the Chernoff function) is defined through 
\begin{equation}
\varphi ^{W}(x)=M^{\ast }(x):=\sup_{t}tx-M(t).  \label{phiWduale de MgfW}
\end{equation}%
The function $x\rightarrow \varphi ^{W}(x)$ is non negative, is $C^{\infty }$
and convex. We also assume that $EW=1$ together with $VarW=1$ which implies $%
\varphi ^{W}(1)=\left( \varphi ^{W}\right) ^{\prime }(1)=0$ and $\left(
\varphi ^{W}\right) ^{\prime \prime }(1)=1.$ Hence $\varphi ^{W}(x)$ is a
divergence function with corresponding divergence pseudo-distance $\phi
^{W}. $ Associated with $\varphi ^{W}$ is the conjugate divergence $%
\widetilde{\phi ^{W}}$ with divergence function $\widetilde{\varphi ^{W}}$ ,
which therefore satisfies $\phi ^{W}\left( Q,P\right) =\widetilde{\phi ^{W}}%
\left( P,Q\right) $ whenever neither $P$ nor $Q$ have null entries.

It is \ of interest to note that the classical power divergences $%
\varphi_{\gamma}$ can be represented through (\ref{phiWduale de MgfW}) for $%
\gamma\leq1$ or $\gamma\geq2.$ A first proof of this lays in the fact that
when $W$ has a distribution in a Natural Exponential Family (NEF) with power
variance function $\alpha=2-\gamma$, then the Legendre transform $%
\varphi^{W} $ of its cumulant generating function $M$ is indeed of the form (%
\ref{CRdiv}). See \cite{LetacMora90} and \cite{BarLevEnis86} for NEF's and
power variance functions, and \cite{Br2017} for relation to the bootstrap.\
A general result of a different nature , including the former ones, can be
seen in \cite{BrSt2020}, Theorem 20. Correspondence between the various
values of $\gamma$ and the distribution of the respective weights can be
found in \cite{BrSt2020}, Example 39, and it can be summarized as presented
now.

For $\gamma <0$ the RV $W$ is constructed as follows: Let $Z$ be an
auxiliary RV with density $f_{Z\text{ }}$ and support $[0,\infty )$ of a
stable law with parameter triplet $\left( -\frac{\gamma }{1-\gamma },0,\frac{%
(1-\gamma )^{-\gamma //(1-\gamma )}}{\gamma }\right) $ in terms of the "form
B notation" on p 12 in \cite{Zolo}, and 
\begin{equation*}
f_{W}(y):=\frac{\exp \left( -y/(1-\gamma \right) )}{\exp (1/\gamma )}f_{Z}(y)%
{\large 1}_{[0,\infty )}(y).
\end{equation*}%
For $\gamma =0$ (which amounts to consider the limit as $\gamma \rightarrow
0 $ in (\ref{CRdiv})) then $W$ has a standard exponential distribution $E(1)$
on $[0,\infty ).$

For $\gamma\in\left( 0,1\right) $ then $W$ has a compound Gamma-Poisson
distribution $C\left( POI(\theta),GAM(\alpha,\beta)\right) $ where $%
\theta=1/\gamma$, $\alpha=1/(1-\gamma)$ and $\beta=\gamma/(1-\gamma).$

For $\gamma=1$ then $W$ has a Poisson distribution with parameter $1,POI(1).$

For $\gamma=2$ then $W$ has normal distribution with expectation and
variance equal to $1.$

For $\gamma >2$ then the RV $W$ is constructed as follows: Let $Z$ be an
auxiliary RV with density $f_{Z\text{ }}$ and support $(-\infty ,\infty )$
of a stable law with parameter triplet $\left( \frac{\gamma }{\gamma -1},0,%
\frac{(\gamma -1)^{-\gamma //(\gamma -1)}}{\gamma }\right) $ in terms of the
"form B notation" on p 12 in \cite{Zolo}, and 
\begin{equation*}
f_{W}(y):=\frac{\exp \left( y/(\gamma -1\right) )}{\exp (1/\gamma )}f_{Z}(-y)%
\text{{\large \ \ \ \ ,}}y\in \mathbb{R}.
\end{equation*}

\bigskip

\section{Maximum likelihood under finitely supported distributions and
simple sampling\label{Sect ML under finite supported}}

\subsection{Standard derivation\label{Subsect Standard derivatioin}}

Let $X_{1},...X_{n}$ be a set of $n$ independent random variables with
common probability measure $P_{\theta _{T}}$ and consider the Maximum
Likelihood estimator of $\theta _{T}$ . A common way to define the ML
paradigm is as follows: For any $\theta $ consider independent random
variables $\left( X_{1,\theta },...X_{n,\theta }\right) $ with probability
measure $P_{\theta }$ , thus \textit{sampled in the same way as the }$X_{i}$%
\textit{'s}, but under some alternative $\theta .$ Define $\theta _{ML}$ as
the value of the parameter $\theta $ for which the probability that, up to a
permutation of the order of the $X_{i,\theta }$'s, the probability that $%
\left( X_{1,\theta },...X_{n,\theta }\right) $ coincides with $%
X_{1},...X_{n} $ is maximal, conditionally on the observed sample $%
X_{1},...X_{n}.$ In formula, let $\sigma $ denote a random permutation of
the indexes $\left\{ 1,2,...,n\right\} $ and $\theta _{ML}$ is defined
through 
\begin{equation}
\theta _{ML}:=\arg \max_{\theta }\frac{1}{n!}\sum_{\sigma \in \mathfrak{S}%
}P_{\theta }\left( \left. \left( X_{\sigma (1),\theta },...,X_{\sigma
(n),\theta }\right) =\left( X_{1},...X_{n}\right) \right\vert \left(
X_{1},...X_{n}\right) \right)  \label{MLfinite}
\end{equation}%
where the summation is extended on all equally probable permutations of $%
\left\{ 1,2,...,n\right\} .$

Denote

\begin{equation*}
P_{n}:=\frac{1}{n}\sum_{i=1}^{n}\delta_{X_{i}}
\end{equation*}
and 
\begin{equation*}
P_{n,\theta}:=\frac{1}{n}\sum_{i=1}^{n}\delta_{X_{i,\theta}}
\end{equation*}
the empirical measures pertaining respectively to $\left(
X_{1},...X_{n}\right) $ and $\left( X_{1,\theta},...X_{n,\theta}\right) $

An alternative expression for $\theta_{ML}$ is 
\begin{equation}
\theta_{ML}:=\arg\max_{\theta}P_{\theta}\left( \left.
P_{n,\theta}=P_{n}\right\vert P_{n}\right) .  \label{MLfinitebis}
\end{equation}

\bigskip An explicit enumeration of the above expression $P_{\theta }\left(
\left. P_{n,\theta }=P_{n}\right\vert P_{n}\right) $ involves the quantities 
\begin{equation*}
n_{j}:=card\left\{ i:X_{i}=d_{j}\right\} 
\end{equation*}%
for $j=1,...,K$ and yields%
\begin{equation}
P_{\theta }\left( \left. P_{n,\theta }=P_{n,X}\right\vert P_{n,X}\right) =%
\frac{n!P_{\theta }\left( d_{j}\right) ^{n_{j}}}{\prod%
\limits_{j=1}^{K}n_{j}!}  \label{multinomial}
\end{equation}%
as follows from the classical multinomial distribution. Optimizing on $%
\theta $ in (\ref{multinomial}) yields 
\begin{align*}
\theta _{ML}& =\arg \max_{\theta }\sum_{j=1}^{K}\frac{n_{j}}{n}\log
P_{\theta }\left( d_{j}\right)  \\
& =\arg \max_{\theta }\frac{1}{n}\sum_{i=1}^{n}\log P_{\theta }\left(
X_{i}\right) .
\end{align*}%
Consider now the Kullback-Leibler distance between $P_{\theta }$ and $P_{n}$
which is non commutative and defined through 
\begin{align}
KL\left( P_{n},\theta \right) & :=\sum_{j=1}^{k}\varphi \left( \frac{n_{j}/n%
}{P_{\theta }\left( d_{j}\right) }\right) P_{\theta }\left( d_{j}\right)  
\notag \\
& =\sum_{j=1}^{k}\left( n_{j}/n\right) \log \frac{n_{j}/n}{P_{\theta }\left(
d_{j}\right) }  \label{KLempfini}
\end{align}%
where 
\begin{equation}
\varphi _{1}(x):=x\log x-x+1  \label{divKL}
\end{equation}%
which is the Kullback-Leibler divergence function. Minimizing the
Kullback-Leibler distance $KL\left( P_{n},\theta \right) $ upon $\theta $
yields 
\begin{align*}
\theta _{KL}& =\arg \min_{\theta }KL\left( P_{n},\theta \right)  \\
& =\arg \min_{\theta }-\sum_{j=1}^{K}\frac{n_{j}}{n}\log P_{\theta }\left(
d_{j}\right)  \\
& =\arg \max_{\theta }\sum_{j=1}^{K}\frac{n_{j}}{n}\log P_{\theta }\left(
d_{j}\right)  \\
& =\theta _{ML}.
\end{align*}%
Introduce the \textit{conjugate divergence function }$\widetilde{\varphi }%
=\varphi _{0}$ $\ $of $\varphi _{1}$ , inducing the modified
Kullback-Leibler, or so-called Likelihood divergence pseudo-distance $KL_{m}$
which therefore satisfies 
\begin{equation*}
KL_{m}\left( \theta ,P_{n}\right) =KL\left( P_{n},\theta \right) .
\end{equation*}%
We have seen that minimizing the Kullback-Leibler divergence $KL\left(
P_{n},\theta \right) $ amounts to minimizing the Likelihood divergence $%
KL_{m}\left( \theta ,P_{n}\right) $ and produces the ML estimate of $\theta
_{T}.$

\subsection{Asymptotic derivation \label{Sub Asymptotic derivation}}

We assume that 
\begin{equation*}
\lim_{n\rightarrow \infty }P_{n}=P_{\theta _{T}}\text{ \ \ \ a.s.}
\end{equation*}%
This holds for example when the $X_{i}$'s are drawn as an iid sample with
common law $P_{\theta _{T}}$ which we may assume in the present context$.$
From an asymptotic standpoint, Kullback-Leibler divergence is related to the
way $P_{n}$ keeps away from $P_{\theta }$ when $\theta $ is not equal to the
true value of the parameter $\theta _{T}$ generating the observations $X_{i}$%
's and is closely related with the type of sampling of the $X_{i}$'s. In the
present case, when i.i.d. sampling of the $X_{i,\theta }$'s under $P_{\theta
}$ are performed, Sanov Large Deviation theorem leads to 
\begin{equation}
\lim_{n\rightarrow \infty }\frac{1}{n}\log P_{\theta }\left( \left.
P_{n,\theta }=P_{n}\right\vert P_{n}\right) =-KL\left( \theta _{T},\theta
\right) .  \label{Sanov fini}
\end{equation}%
This result can easily be obtained from (\ref{multinomial}) using Stirling
formula to handle the factorial terms and the law of large numbers which
states that for all $j$'s, $n_{j}/n$ tends to $P_{\theta _{T}}(d_{j})$ as $n$
tends to infinity. Comparing with (\ref{KLempfini}) we note that the MLE $%
\theta _{ML}$ is a proxy of the minimizer of the natural estimator $\theta
_{T}$ of $KL\left( \theta _{T},\theta \right) $ in $\theta ,$ substituting
the unknown measure generating the $X_{i}$'s by its empirical counterpart $%
P_{n}$ . Alternatively as will be used in the sequel, $\theta _{ML}$
minimizes upon $\theta $ the Likelihood divergence $KL_{m}\left( \theta
,\theta _{T}\right) $ between $P_{\theta }$ and $P_{\theta _{T}}$
substituting the unknown measure $P_{\theta _{T}}$ generating the $X_{i}$'s
by its empirical counterpart $P_{n}$ . Summarizing we have obtained:

The ML\ estimate can be obtained from a LDP\ statement as given in (\ref%
{Sanov fini}), optimizing in $\theta$ in the estimator of the LDP rate where
the plug-in method of the empirical measure of the data is used instead of
the unknown measure $P_{\theta_{T}}.$ Alternatively it holds 
\begin{equation}
\theta_{ML}:=\arg\min_{\theta}\widehat{KL_{m}}\left( \theta,\theta
_{T}\right)  \label{ML finite case}
\end{equation}
with 
\begin{equation*}
\widehat{KL_{m}}\left( \theta,\theta_{T}\right) :=KL_{m}\left( \theta
,P_{n}\right) .
\end{equation*}

This principle will be kept throughout this paper: the estimator is defined
as maximizing the probability that the simulated empirical measure be close
to the empirical measure as observed on the sample, conditionally on it,
following the same sampling scheme. This yields a maximum likelihood
estimator, and its properties are then obtained when randomness is
introduced as resulting from the sampling scheme.

\section{Bootstrap and weighted sampling\label{Sect Bootstrap and weighted
sampling}}

The sampling scheme which we consider is commonly used in connection with
the bootstrap and is referred to as the \textit{weighted} or \textit{%
generalized bootstra}p, sometimes called \textit{wild bootstrap}, first
introduced by Newton and Mason \cite{Mason}.

Let $X_{1},...,X_{n}$ with common distribution $P_{\theta_{T}}$ on $\mathcal{%
Y}:=\left\{ d_{1},..,d_{K}\right\} .$

Consider a collection $W_{1},...,W_{n}$ of independent copies of $W$, whose
distribution satisfies the conditions stated in Section 1. The weighted
empirical measure $P_{n}^{W}$ is defined through

\begin{equation*}
P_{n}^{W}:=\frac{1}{n}\sum_{i=1}^{n}W_{i}\delta _{X_{i}}.
\end{equation*}%
This empirical measure need not be a probability measure, since its mass may
not equal $1.$ Also it might not be positive, since the weights may take
negative values.\ Therefore $P_{n}^{W}$ can be identified with a random
point in $\mathbb{R}^{K}.$ The measure $P_{n}^{W}$ converges almost surely
to $P_{\theta _{T}}$ when the weights $W_{i}$'s satisfy the hypotheses
stated in Section 1. Indeed general results pertaining to this sampling
procedure state that under regularity, functionals of the measure $P_{n}^{W}$
are asymptotically distributed as are the same functionals of $P_{n}$ when
the $X_{i}$'s are i.i.d. Therefore the weighted sampling procedure mimics
the i.i.d. sampling fluctuation in a two steps procedure: choose $n$ values
of $X_{i}$ such that they asymptotically fit to $P_{\theta _{T}}$, which
means%
\begin{equation*}
\lim_{n\rightarrow \infty }\frac{1}{n}\sum_{i=1}^{n}\delta
_{X_{i}}=P_{\theta _{T}}
\end{equation*}%
a.s. and then play the $W_{i}$'s on each of the $x_{i}$'s. Then get $%
P_{n}^{W}$, a proxy to the random empirical measure $P_{n}$ .

We also consider the normalized weighted empirical measure 
\begin{equation}
\mathfrak{P}_{n}^{W}:=\sum_{i=1}^{n}Z_{i}\delta _{X_{i}}  \label{P_n^W}
\end{equation}%
where 
\begin{equation}
Z_{i}:=\frac{W_{i}}{\sum_{j=1}^{n}W_{j}}  \label{Z_i}
\end{equation}%
whenever $\sum_{j=1}^{n}W_{j}\neq 0$, and 
\begin{equation*}
\mathfrak{P}_{n}^{W}=\infty
\end{equation*}%
when $\sum_{j=1}^{n}W_{j}=0$, where $\mathfrak{P}_{n}^{W}=\infty $ means $%
\mathfrak{P}_{n}^{W}(d_{k})=\infty $ for all $d_{k}$ in $\mathcal{Y}.$

\bigskip

\subsection{A conditional Sanov type result for the weighted empirical
measure\label{Subsect Conditional Sanov weighted emp meas}}

\bigskip We now state a conditional Sanov type result for the family of
random measures $\mathfrak{P}_{n}^{W}.$ It follows readily from a companion
result pertaining to $P_{n}^{W}$ and enjoys a simple form when the weights $%
W_{i}$ are associated to power divergences, as defined in Section \ref%
{Subsection Weights}. We quote the following results, referring to \cite%
{BrSt2020}.

Consider a set $\Omega$ in $\mathbb{R}^{K}$ such that 
\begin{equation}
cl\Omega=cl(Int\Omega)  \label{Omega regul}
\end{equation}
which amounts to a regularity assumption (obviously met when $\Omega$ is an
open set), which allows for the replacement of the usual $\lim\inf$ and $%
\lim\sup$ by standard limits in usual LDP statements. We denote by $P^{W}$
the probability measure of the random family of iid weights $W_{i}.$

It then holds

\begin{proposition}
\label{THM9}(Theorem 9 in \cite{BrSt2020})The weighted empirical measure $%
P_{n}^{W}$ satisfies a conditional Large Deviation Principle in $\mathbb{R}%
^{K}$ namely, denoting $P$ the a.s. limit of $P_{n},$ 
\[
\lim_{n\rightarrow\infty}\frac{1}{n}\log P^{W}\left( \left.
P_{n}^{W}\in\Omega\right\vert X_{1}^{n}\right) =-\phi^{W}\left(
\Omega,P\right)
\]
where $\phi^{W}\left( \Omega,P\right) :=\inf_{Q\in\Omega}\phi^{W}\left(
Q,P\right) .$
\end{proposition}

As a direct consequence of the former result, it holds, for any $\Omega
\subset\mathbb{S}^{K}$ satisfying (\ref{Omega regul}), where $\mathbb{S}^{K}$
designates the simplex of all pm's on $\mathcal{Y}$

\begin{theorem}
\label{THM12}(Theorem 12 in \cite{BrSt2020})The normalized weighted
empirical measure $\mathfrak{P}_{n}^{W}$ satisfies \ a conditional Large
Deviation Principle in $\mathbb{S}^{K}$ 
\begin{equation}
\lim_{n\rightarrow\infty}\frac{1}{n}\log P^{W}\left( \left. \mathfrak{P}%
_{n}^{W}\in\Omega\right\vert X_{1}^{n}\right) =-\inf_{m\neq0}\phi^{W}\left(
m\Omega,P\right) .  \label{Sanov weighted  normalized}
\end{equation}
\end{theorem}

A flavour of the simple proofs of Proposition \ref{THM9} and Theorem \ref%
{THM12} is presented in the Appendix; see \cite{BrSt2020} for a detailed
treatment; see also Theorem 3.2\ and Corollary 3.3 in \cite{WiTr2014} where
Theorem \ref{THM12} is proved in a more abstract setting. Note that the
mapping $Q\rightarrow \inf_{m\neq 0}\phi ^{W}\left( mQ,P\right) $ is indeed
a divergence in the simplex $\mathbb{S}^{K}$ for all pm $P$ defined on $%
\mathcal{Y}$ with positive entries.\ 

We will be interested in the pm's in $\Omega$ which minimize the RHS in the
above display. The case when $\ \phi^{W}$ is a power divergence, namely $%
\phi^{W}=\phi_{\gamma}$ for some real $\gamma\in(-\infty,1]\cup
\lbrack2,\infty)$ enjoys a special property with respect to the pm's $Q$
achieving the infimum (upon $Q$ in $\Omega$) in (\ref{Sanov weighted
normalized}). It holds

\begin{proposition}
\label{PropoLemma14}(Lemma 14 in \cite{BrSt2020})Assume that $\phi^{W}$ is a
power divergence. Then 
\[
Q\in\arg\inf\left\{ \inf_{m\neq0}\phi^{W}\left( mQ,P\right) ,Q\in
\Omega\right\}
\]
and 
\[
Q\in\arg\inf\left\{ \phi^{W}\left( Q,P\right) ,Q\in\Omega\right\}
\]
are equivalent statements.
\end{proposition}

Indeed Proposition \ref{PropoLemma14} holds as a consequence of the
following results, to be used later on.

\begin{lemma}
\label{Lemma rates inf}For $Q$ and $P$ two pm's such that the involved
expressions are finite, it holds

(i) For $\gamma\in\left( 0,1\right) $ it holds $\inf_{m\neq0}\phi_{\gamma
}(mQ,P)=\left( 1-\gamma\right) \phi_{\gamma}\left( Q,P\right) .$

(ii) For $\gamma<0$ and $\gamma>1$ it holds $\inf_{m\neq0}\phi_{\gamma
}(mQ,P)=\frac{1}{\gamma}\left[ 1-\left( 1+\gamma(\gamma-1)\phi_{\gamma
}(Q,P)\right) ^{-1/(\gamma-1)}\right] .$

(iii) $\inf_{m\neq0}\phi_{1}(mQ,P)=1-\exp\left( -KL(Q,P)\right)
=1-\exp(-\phi_{1}(Q,P)).$

(iv)$\inf_{m\neq0}\phi_{0}(mQ,P)=KL_{m}(Q,P)=\phi_{0}(Q,P)$
\end{lemma}

The weighted empirical measure $P_{n}^{W}$ has been used in the weighted
bootstrap (or wild bootstrap) context, although it is not a pm.\ However,
conditionally upon the sample points, its produces statistical estimators $%
T(P_{n}^{W})$ whose weak behavior (conditionally upon the sample) converges
to the same limit as does $T(P_{n})$ when normalized on the classical CLT
range; see eg Newton and Mason \cite{Mason}. Large deviation theorem for the
weighted empirical measure $P_{n}^{W}$ has been obtained by \cite{Barbe95};
for other contributions in line with those, see \cite{Najim} and \cite%
{WiTr2014}.\ Normalizing the weights produces families of exchangeable
weights $Z_{i}$, and the normalized weighted empirical measure $\mathfrak{P}%
_{n}^{W}$ is the cornerstone for the so-called non parametric Bayesian
bootstrap, initiated by \cite{Rubin81}, and further developed by \cite%
{NewRaft} among others. Note however that in this context the RV's $W_{i}$'s
are chosen as distributed as standard exponential variables.\ The link with
spacings from a uniform distribution and the corresponding reproducibility
of the Dirichlet distributions are the basic ingredients which justify the
non parametric bootstrap approach; in the present context, the choice of the
distribution of the $W_{i}$'s is a natural extension of this paradigm, at
least when those $W_{i}$'s are positive RV's.\bigskip

\bigskip

\subsection{Maximum Likelihood for the generalized bootstrap\label{Subsect
ML generalized bootstrap}}

We will consider maximum likelihood in the same spirit as developed in
Section \ref{Sub Asymptotic derivation}, here in the context of the
normalized weighted empirical measure; it amounts to justify minimum
divergence estimators as appropriate MLE's under such bootstrap procedure.\ 

We thus consider the same statistical model $\mathcal{P}_{\Theta }$ and keep
in mind the ML principle as seen as resulting from a maximization of the
conditional probability of getting simulated observations close to the
initially observed data. Similarly as in Section \ref{Sect ML under finite
supported} fix an arbitrary $\theta $ and simulate $X_{1,\theta
},..,X_{n,\theta }$ with distribution $P_{\theta }.$ Define accordingly $%
P_{n,\theta }^{W}$ and $\mathfrak{P}_{n,\theta }^{W}$ making use of iid RV's 
$W_{1},..,W_{n}$ . Now the event $\mathfrak{P}_{n,\theta }^{W}(k)=n_{k}/n$
has probability $0$ in most cases (for example when $W$ has a continuous
distribution), and therefore we are led to consider events of the form $%
\mathfrak{P}_{n,\theta }^{W}\in V_{\varepsilon }\left( P_{n}\right) $ ,
meaning $\max_{k}\left\vert \mathfrak{P}_{n,\theta
}^{W}(d_{k})-P_{n}(d_{k})\right\vert \leq \varepsilon $ for some positive $%
\varepsilon ;$ notice that $V_{\varepsilon }\left( P_{n}\right) $ defined
through $V_{\varepsilon }\left( P_{n}\right) :=\left\{ Q\in \mathbb{S}%
^{K}:\max_{k}\left\vert Q(d_{k})-P_{n}(d_{k})\right\vert \leq \varepsilon
\right\} $ has non void interior$.$

For such a configuration consider 
\begin{equation}
P^{W}\left( \left. \mathfrak{P}_{n,\theta}^{w}\in V_{\varepsilon}\left(
P_{n}\right) \right\vert X_{1,\theta},..,X_{n,\theta},X_{1},..,X_{n}\right)
\label{ML Bootstrap}
\end{equation}
where the $X_{i,\theta}$ are randomly drawn iid under $P_{\theta}.$
Obviously for $\theta$ far away from $\theta_{T}$ the sample $\left(
X_{1,\theta },..,X_{n,\theta}\right) $ is realized "far away " from $\left(
X_{1},..,X_{n}\right) $, which has been generated under the truth, namely $%
P_{\theta_{T}}$ , and the probability in (\ref{ML Bootstrap}) is small,
whatever the weights, for small $\varepsilon$.

\bigskip

We will now consider (\ref{ML Bootstrap}) asymptotically on $n$, since, in
contrast with the first derivation of the standard MLE in Section \ref%
{Subsect Standard derivatioin}, we cannot perform the same calculation for
each $n$, which was based on multinomial counts. Note that we obtained a
justification for the usual MLE through the asymptotic Sanov LDP, leading to
the KL divergence and finally back to the MLE through an approximation step
of this latest.

\bigskip

We first state

\begin{theorem}
With the above notation the following conditioned LDP result holds, for some 
$\alpha<1<\beta$

\begin{align}
-\inf_{m\neq0}\phi^{W}(mV_{\alpha\epsilon}(P_{\theta_{T}}),\theta) &
\leq\lim_{n\rightarrow\infty}\frac{1}{n}\log P^{W}\left( \mathfrak{P}%
_{n,\theta}^{W}\in
V_{\epsilon}(P_{n})|X_{1,\theta},...,X_{n,\theta},X_{1},..,X_{n}\right)
\label{LDP voisinage P_n} \\
& \leq-\inf_{m\neq0}\phi^{W}(mV_{\beta\epsilon}(P_{\theta_{T}}),\theta ) 
\nonumber
\end{align}
where $\phi^{W}(V_{c\epsilon}(\theta_{T}),\theta)=\inf_{\mu\in V_{c\epsilon
}(P_{\theta_{T}}))}\phi^{W}(\mu,\theta)$.
\end{theorem}

\bigskip The above result follows from Theorem \ref{Sanov weighted
normalized} together with the a.s. convergence of $P_{n}$ to $P_{\theta_{T}}$
in $\mathbb{S}^{K}.$

From the above result it appears that as $\varepsilon\rightarrow0$ , by
continuity it holds%
\begin{equation}
\lim_{\varepsilon\rightarrow0}\lim_{n\rightarrow\infty}\frac{1}{n}\log
P^{W}\left( \mathfrak{P}_{n,\theta}^{W}\in V_{\epsilon}(P_{n})|X_{1,\theta
},...,X_{n,\theta},X_{1},..,X_{n}\right)
=-\inf_{m\neq0}\phi^{W}(mP_{\theta_{T}},\theta).
\label{passage limite epsilon}
\end{equation}
The ML principle amounts to maximize $P^{W}\left( \mathfrak{P}%
_{n,\theta}^{W}\in
V_{\epsilon}(P_{n})|X_{1,\theta},...,X_{n,\theta},X_{1},..,X_{n}\right) $
upon $\theta.$ Whenever $\Theta$ is a compact set we may insert this
optimization in (\ref{LDP voisinage P_n}) which yields, following (\ref%
{passage limite epsilon})%
\begin{equation*}
\lim_{\varepsilon\rightarrow0}\lim_{n\rightarrow\infty}\frac{1}{n}\log
\sup_{\theta}P^{W}\left( \mathfrak{P}_{n,\theta}^{W}\in
V_{\epsilon}(P_{n})|X_{1,\theta},...,X_{n,\theta},X_{1},..,X_{n}\right)
=-\inf_{\theta \in\Theta}\inf_{m\neq0}\phi^{W}(mP_{\theta_{T}},\theta).
\end{equation*}

By Proposition \ref{PropoLemma14} the argument of the infimum upon $\theta $
in the RHS\ of the above display coincides with the corresponding argument
of $\phi ^{W}(\theta _{T},\theta )$, which obviously gets $\theta _{T}$.
This justifies to consider a proxy of this minimization problem as a "ML"
estimator based on normalized weighted data.\ 

A further interpretation of the MDE as a Maximum a posteriory estimator
(MAP) \ in the context of non parametric bayesian procedures may also be
proposed; this is postponed to a next paper.

Since 
\begin{equation*}
\phi ^{W}(\theta _{T},\theta )=\widetilde{\phi }^{W}(\theta ,\theta _{T})
\end{equation*}%
the ML\ estimator is obtained as in the conventional case by plug in the LDP
rate.\ Obviously the "best" plug in consists in the substitution of $%
P_{\theta _{T}}$ by $P_{n},$ the empirical measure of the sample, since $%
P_{n}$ achieves the best rate of convergence to $P_{\theta _{T}}$ when
confronted to any bootstrapped version, which adds "noise" to the sampling.\
We may therefore call 
\begin{align}
\theta _{ML}^{W}& :=\arg \inf_{\theta \in \Theta }\widetilde{\phi }%
^{W}(\theta ,P_{n}):=\arg \inf_{\theta \in \Theta }\sum_{k=1}^{K}P_{n}(d_{k})%
\widetilde{\varphi }\left( \frac{P_{\theta }(d_{k})}{P_{n}(d_{k})}\right)
\label{def avec P_n} \\
& =\arg \inf_{\theta \in \Theta }\sum_{k=1}^{K}P_{\theta }(d_{k})\varphi
\left( \frac{P_{n}(d_{k})}{P_{\theta }(d_{k})}\right)  \notag
\end{align}%
the MLE\ for the bootstrap sampling; here $\widetilde{\phi }^{W}$ (with
divergence function $\widetilde{\varphi }$) is the conjugate divergence of $%
\phi ^{W}$ (with divergence function $\varphi $) .Since $\phi ^{W}=\phi
_{\gamma }$ for some $\gamma $, it holds $\widetilde{\phi }^{W}=\phi
_{1-\gamma }.$

\bigskip Obviously we can also plug in the normalized weighted empirical
measure, which also is a proxy of $P_{\theta _{T}}$ for each run of the
weights.\ This produces a bootstrap estimate of $\theta _{T}$ through 
\begin{align}
\theta _{B}^{W}& :=\arg \inf_{\theta \in \Theta }\widetilde{\phi }%
^{W}(\theta ,\mathfrak{P}_{n}^{W}):=\arg \inf_{\theta \in \Theta
}\sum_{k=1}^{K}\mathfrak{P}_{n}^{W}(d_{k})\widetilde{\varphi }\left( \frac{%
P_{\theta }(d_{k})}{\mathfrak{P}_{n}^{W}(d_{k})}\right)
\label{def avec P_n boot} \\
& =\arg \inf_{\theta \in \Theta }\sum_{k=1}^{K}P_{\theta }(d_{k})\varphi
\left( \frac{\mathfrak{P}_{n}^{W}(d_{k})}{P_{\theta }(d_{k})}\right)  \notag
\end{align}

\bigskip

where $\mathfrak{P}_{n}^{W}$ is defined in (\ref{P_n^W}), assuming $n$ large
enough such that this the sum of the $W_{i}$'s not zero. Whenever $W$ has
positive probability to assume value $0$, these estimators are defined for
large $n$ in order that $\mathfrak{P}_{n}^{W}(d_{k})$ be positive for all $%
k. $ Since $E(W)=1$, this occurs for large samples.

When $\mathcal{Y}$ is not a finite space then an equivalent construction can
be developed based on the variational form of the divergence; see \cite%
{BK2009}.

\begin{remark}
We may also consider cases when the MLE defined through $\theta _{ML}^{W}$
defined in (\ref{def avec P_n}) coincide with the standard MLE $\theta _{ML}$
under iid sampling, and when their bootstrapped counterparts $\theta
_{B}^{W} $ defined in (\ref{def avec P_n boot}) coincides with the
bootstrapped standard MLE $\theta _{ML}^{b}$ defined through the likelihood
estimating equation where the factor $1/n$ is substituted by the weight $%
Z_{i}$ .\ It is proved \ in Theorem 5 of \cite{BR2014} that whenever $%
\mathcal{P}_{\Theta }$ is an exponential family with natural parametrization 
$\theta \in \mathbb{R}^{d}$ and sufficient statistics $T$%
\begin{equation*}
P_{\theta }\left( d_{j}\right) =\exp \left[ T(d_{j})^{\prime }\theta
-C(\theta )\right] ,\ \ \ \ 1\leq j\leq K
\end{equation*}%
where the Hessian matrix of $C(\theta )$ is definite positive, then for all
divergence pseudo distance $\phi $ satisfying \ regularity conditions
(including therefore the present cases), $\theta _{ML}^{W}$ equals $\theta
_{ML}$, the classical MLE\ in $\mathcal{P}_{\Theta }$ defined as the
solution of the normal equation%
\begin{equation*}
\frac{1}{n}\sum T(X_{i})=\nabla C(\theta _{ML})
\end{equation*}%
irrespectively upon $\phi .$ Therefore on regular exponential families, and
under iid sampling, all minimum divergence estimators coincide with the MLE
(which is indeed one of them).\ The proof of this result is based on the
variational form of the estimated of divergence $Q\rightarrow \phi \left(
Q,P\right) $, which coincides with the plug in version in (\ref{def avec P_n}%
) when the common support of all distribution in $\mathcal{P}_{\Theta }$ is
finite.\ \ Following verbatim the proof of Theorem 5 in \cite{BR2014}\
substituting $P_{n}$ by $\mathfrak{P}_{n}^{W}$ it results that $\theta
_{B}^{W}$ equals the weighted MLE (standard generalized bootstrapped MLE $%
\theta _{ML}^{b}$) defined through the normal equation 
\begin{equation*}
\sum_{i=1}^{n}Z_{i}T(X_{i})=\nabla C(\theta _{ML}^{b}).
\end{equation*}%
where the $Z_{i}$'s are defined in (\ref{Z_i}).\ This fact holds for any
choice of the weights, irrespectively on the choice of the divergence
function $\varphi $ with the only restriction that it satisfies the mild
conditions (RC) in \cite{BR2014}. It results that for those models any
generalized bootstrapped MDE coincides with the corresponding bootstrapped
MLE.
\end{remark}

\begin{remark}
The \ estimators $\theta _{ML}^{W}$ defined in (\ref{def avec P_n}) have
been considered for long irrespectively of the present approach; see e.g. 
\cite{Pardo}.\ Their statistical properties in various contexts have been
studied \ for general support spaces $\mathcal{Y}$ in \cite{BK2009} for
parametric models, and for various semi parametric models in \cite{BrK2012}
and \cite{BrDec2016}.
\end{remark}

\begin{example}
\label{Example distr W vs div}

A-In \ the case when $W$ is a RV with standard exponential distribution,
then the normalized weighted empirical measure $\mathfrak{P}_{n}^{W}$ is a
realization of the a posteriori distribution for the non informative prior
on the non parametric distribution of $X.$ See \cite{Rubin81}. In this case $%
\varphi (x)=-\log x+x-1$ and $\widetilde{\varphi }(x)=x\log x-x+1;$ the
resulting estimator is the minimum Kullback-Leibler one.\ 

B-When $W$ has a standard Poisson distribution then the couple $\left(
\varphi ,\widetilde{\varphi }\right) $ is reverse wrt the above one, and the
resulting estimator is the minimum modified Kullback-Leibler one. which
takes the usual weighted form of the standard generalized bootstrap \ MLE%
\begin{equation*}
\theta _{B}^{POI(1)}:=\arg \sup_{\theta }\sum_{k=1}^{K}\left( \frac{%
\sum_{i=1}^{n}W_{i}{\Large 1}_{k}(X_{i})}{\sum_{i=1}^{n}W_{i}}\right) \log
P_{\theta }(k)
\end{equation*}%
which is defined for $n$ large enough so that $\sum_{i=1}^{n}W_{i}\neq 0.$
Also in this case $\theta _{ML}^{W}$ coincides with the standard MLE.

C-In case when $W$ has an Inverse Gaussian distribution IG(1,1) then $\
\varphi (x)=\varphi _{-1}(x)=\frac{1}{2}\left( x-1\right) ^{2}/x$ for $x>0$
and the ML estimator minimizes the Pearson Chi-square divergence with
generator function $\varphi _{2}(x)=\frac{1}{2}\left( x-1\right) ^{2}$ which
is defined on $\mathbb{R}.$

D-When $W$ follows a normal distribution with expectation and variance $1$,
then the resulting divergence is the Pearson Chi-square divergence $\varphi
_{2}(x)$ and the resulting estimator minimizes the Neyman Chi-square
divergence with $\varphi (x)=\varphi _{-1}(x).$

E-When $W$ has a Compound Poisson Gamma distribution $C\left( POI(2),\Gamma
(2,1)\right) $ distribution then the corresponding divergence is $\varphi
_{1/2}(x)=2\left( \sqrt{x}-1\right) ^{2}$ which is self conjugate, whence
the ML\ estimator is the minimum Hellinger distance one.
\end{example}

\section{Optimal weighting in relation with the bootstrapped estimator of
the divergence\label{Sect Optimal weighting}}

The definition of the divergence estimator (\ref{def avec P_n boot}) opens
to the definition of various bootstrap estimators for the divergence pseudo
distance between $P_{\theta _{T}}$ and the model $\mathcal{P}.$ Indeed the
choice of the distribution of the weights $W_{i}$'s in $\mathfrak{P}_{n}^{W}$
needs not be related to the divergence function $\phi $; for example we
might define some $\mathfrak{P}_{n}^{V}$ in order to define an estimator of $%
\phi ^{W}(Q,P)$ through (\ref{def avec P_n boot}) with $\mathfrak{P}_{n}^{W}$
substituted by $\mathfrak{P}_{n}^{V}$ where 
\begin{equation*}
Z_{i}:=\frac{V_{i}}{\sum_{j=1}^{n}V_{j}}
\end{equation*}%
and the vector $\left( V_{1},..,V_{n}\right) $ is not related in any way
with $\phi ^{W}$. The resulting generalized bootstrapped estimate of $\phi
^{W}(Q,P)$ may also be used as a test statistics in order to assess whether $%
\theta =\theta _{T}$ for example.\ It may seem a natural insight that the
choice when $\left( V_{1},..,V_{n}\right) $ has same distribution as $\left(
W_{1},..,W_{n}\right) $ should bear some optimality property.\ The r\^{o}le
of the present section is to explore this question, for specific choices of
\ the distribution of the $V_{i}$'s.

\subsection{Comparing bootstrapped statistics\label{Subsect Comparing
bootstrapped stats}}

Let $\phi ^{W}$ be a power divergence defined by some weight $W$ through (%
\ref{phiWduale de MgfW})$.$ We assume that $\theta _{T}$ is known and we
measure the divergence between $\mathfrak{P}_{n}^{W}$ and $P_{\theta _{T}}$
as a bootstrapped version of the corresponding distance between $P_{n}$ and $%
P_{\theta _{T}}$, where the distance is suited to the distribution of the
weights. We compare the decay to $0$ of this same distance with the
corresponding decay substituting $\mathfrak{P}_{n}^{W}$ by $\mathfrak{P}%
_{n}^{V}$ for some competing family of weights $\left( V_{1},..,V_{n}\right)
.$ Both RV's $W$ and $V$ are assumed to have distributions such that the
Legendre transform of their cumulant generating functions belong to the
Cressie Read family of divergences.\ The divergence $\phi ^{W}$ is
associated to the generator $\varphi _{\gamma }$ and, respectively, $V$ is
associated to a generator $\varphi _{\gamma ^{\prime }}$ by the
corresponding formula (\ref{phiWduale de MgfW})$.$ \bigskip For brevity we
restrict the discussion to RV's $W$ and $V$ which are associated to
divergence functions $\varphi _{\gamma }$ and $\varphi _{\gamma ^{\prime }}$
with $\gamma $ $\in \left( 0,1\right) $, as other cases are similar, making
use of the corresponding formulas from Lemma \ref{Lemma rates inf}.

The bootstrap distance between $P_{\theta_{T}}$ and the bootstrapped dataset
will be defined $\ \ $as$\ \ \ \widetilde{\phi_{\gamma}}\left( \theta _{T},%
\mathfrak{P}_{n}^{W}\right) .$

Looking at case (i) in Proposition \ref{PropoLemma14} we denote by $\varphi$
the generator of the divergence $Q\rightarrow\inf_{m\neq0}\phi_{%
\gamma}(mQ,P) $ and $\psi$ the generator of the divergence $%
Q\rightarrow\inf_{m\neq0}\phi_{\gamma^{\prime}}(mQ,P)$ from which 
\begin{equation*}
\varphi(x)=\frac{x^{\gamma}-\gamma x+\gamma-1}{\gamma}
\end{equation*}%
\begin{equation*}
\psi(x)=\frac{x^{\gamma^{\prime}}-\gamma^{\prime}x+\gamma^{\prime}-1}{%
\gamma^{\prime}}.
\end{equation*}

Also for clearness we denote 
\begin{equation*}
\Phi(Q,P):=\inf_{m\neq0}\phi_{\gamma}(mQ,P)
\end{equation*}
and 
\begin{equation*}
\Psi\left( Q,P\right) :=\inf_{m\neq0}\phi_{\gamma^{\prime}}(mQ,P)
\end{equation*}
For $\gamma\in\left( 0,1\right) $ due to Proposition \ref{PropoLemma14} (i)
and Lemma \ref{Lemma rates inf}, for $Q$ and $P$ in $\mathbb{S}^{K}$ $\ $\
with non null entries,\ we define the conjugate divergence $\widetilde{\Phi}%
(Q,P):=\Phi(P,Q)$ and the generator of $Q\rightarrow \widetilde{\Phi}(Q,P)$
writes 
\begin{equation*}
\widetilde{\varphi}(x):=\left( \gamma-1\right) \varphi_{1-\gamma}(x)
\end{equation*}
we will denote accordingly $\widetilde{\psi}(x):=\left( \gamma^{\prime
}-1\right) \varphi_{1-\gamma^{\prime}}(x)$ \ the generator of $\widetilde{%
\Psi}(Q,P)$ , the conjugate divergence of $\Psi\left( Q,P\right) .$

\bigskip In order to simplify the notation , for any event $A$, $%
P_{X_{1}^{n}}^{W}$ $(A)$ denotes the probability of $A$ conditioned upon $%
\left( X_{1},..,X_{n}\right) .$\bigskip

\ By Theorem \ref{THM12}%
\begin{equation}
\lim_{n\rightarrow\infty}\frac{1}{n}\log P_{X_{1}^{n}}^{W}\left( \widetilde{%
\phi_{\gamma}}\left( \theta_{T},\mathfrak{P}_{n}^{V}\right) >t\right)
=-\inf\left\{ \Psi\left( Q,\theta_{T}\right) ,Q:\widetilde{\phi _{\gamma}}%
\left( \theta_{T},Q\right) >t\right\} .  \label{LDPQ_n}
\end{equation}

We prove

\begin{proposition}
\label{Prop Bootstrap optim}\ For $\gamma\in\left( 0,1\right) $ and any $%
\gamma^{\prime}\notin\left( 1,2\right) $, for any positive $t$, 
\begin{equation}
\lim_{n\rightarrow\infty}\frac{1}{n}\log P_{X_{1}^{n}}^{W}\left( \widetilde{%
\phi_{\gamma}}\left( \theta_{T},\mathfrak{P}_{n}^{W}\right) >t\right)
=-t(1-\gamma)  \tag{(i)}
\end{equation}
while 
\begin{equation}
\lim_{n\rightarrow\infty}\frac{1}{n}\log P_{X_{1}^{n}}^{W}\left( \widetilde{%
\phi_{\gamma}}\left( \theta_{T},\mathfrak{P}_{n}^{V}\right) >t\right)
\geq-t(1-\gamma).  \tag{(ii)}
\end{equation}
\ 
\end{proposition}

Proof: By Theorem \ref{THM12}, since $\Phi\left( Q,\theta)=(1-\gamma
)\phi_{\gamma}(Q,\theta\right) $, it holds for any $t>0$ 
\begin{align*}
& \lim_{n\rightarrow\infty}\frac{1}{n}\log P_{X_{1}^{n}}^{W}\left( 
\widetilde{\phi_{\gamma}}\left( \theta_{T},\mathfrak{P}_{n}^{W}\right)
>t\right) \\
& =-\inf\left\{ \Phi\left( Q,\theta_{T}\right) ,Q:\widetilde{\phi_{\gamma }}%
\left( \theta_{T},Q\right) >t\right\} \\
& =-\inf\left\{ \Phi\left( Q,\theta_{T}\right) ,Q:\phi_{\gamma}\left(
Q,\theta_{T}\right) >t\right\} \\
& =-\inf\left\{ \Phi\left( Q,\theta_{T}\right) ,Q:\Phi\left( Q,\theta
_{T}\right) >t(1-\gamma)\right\} \\
& =-t(1-\gamma)
\end{align*}
which proves (i).

Using (\ref{LDPQ_n}), for any $t>0$ 
\begin{align*}
\lim_{n\rightarrow\infty}\frac{1}{n}\log P_{X_{1}^{n}}^{W}\left( \widetilde{%
\phi_{\gamma}}\left( \theta_{T},\mathfrak{P}_{n}^{V}\right) >t\right) &
=-\inf\left\{ \Psi\left( Q,\theta_{T}\right) ,Q:\phi _{\gamma}\left(
Q,\theta_{T}\right) >t\right\} \\
& =-\inf\left\{ \Psi\left( Q,\theta_{T}\right) ,Q:\Phi\left( Q,\theta
_{T}\right) >t(1-\gamma)\right\}
\end{align*}
from which (ii) holds whenever there exists some $R$ in $\mathbb{S}^{K}$
satisfying both 
\begin{equation*}
\Phi\left( R,\theta_{T}\right) >t(1-\gamma)
\end{equation*}
and 
\begin{equation*}
\Psi\left( R,\theta_{T}\right) \leq t(1-\gamma).
\end{equation*}

\bigskip

For $K\geq2$, let $R:=\left( r_{1},..,r_{K}\right) $ such that $r_{i}=ap_{i}$
for $i=1,..,K-1$ where $P_{\theta_{T}}:=\left( p_{1},..,p_{K}\right) \in%
\mathbb{S}^{K}$ with non null entries. Assume $a<1.$ Then $%
a\rightarrow\Phi\left( R,P\right) $ is decreasing on $\left( 0,1\right) $, $%
\lim_{a\rightarrow0}$ $\Phi\left( R,P\right) =+\infty$ and $%
\lim_{a\rightarrow1}\Phi\left( R,P\right) =0$; thus there exists $%
a_{\varphi}(t)$ such that for $a\in\left( 0,a_{\varphi}(t)\right) ,$ it
holds $\Phi\left( R,P\right) >t(1-\gamma).$ In the same way there exists $%
a_{\psi}(t)$ such that for $a\in\left( a_{\psi}(t),1\right) $ it holds $%
\Psi\left( R,P\right) <t(1-\gamma).$ Hence for $a\in\left( \min\left(
a_{\varphi}(t),a_{\psi}(t)\right) ,\max\left( a_{\varphi}(t),a_{\psi
}(t)\right) \right) $ , there exists some $R$ which satisfies the claim.

\bigskip

To summarize the meaning of Proposition \ref{Prop Bootstrap optim}, one can
say that it inlights the necessary fit between the divergence and the law of
the weights when exploring the asymptotic behavior of the bootstrapped
empirical measure. It can also be captured stating that given a divergence $%
\phi_{\gamma}$ there exists an optimal bootstrap in the sense that the
chances for $\ \widetilde{\phi_{\gamma}}\left( \theta_{T},\mathfrak{P}%
_{n}^{W}\right) $ to be large are minimal; the "noise" caused by the weights
is tampered down when those are fitted to the divergence, hence in no way in
an arbitrary way.

\subsection{Bahadur efficiency of minimum divergence tests under generalized
bootstrap\label{Subsect Bahadur efficiency}}

In \cite{EfronTibshi} Efron and Tibshirani suggest the bootstrap as a
valuable approach for testing, based on bootstrapped samples.\ We show that
bootstrap testing for parametric models based on appropriate divergence
statistics enjoys maximal Bahadur efficiency with respect to any bootstrap
test statistics.

The standard approach to Bahadur efficiency can be adapted for the present
generalized Bootstrapped tests as follows.

Consider the test of some null hypothesis H0: $\theta_{T}=\theta$ versus a
simple hypothesis H1 $\theta_{T}=\theta^{\prime}.$

We consider two competitive statistics for this problem. The first one is
based on the bootstrap estimate of $\widetilde{\phi }^{W}\left( \theta
,\theta _{T}\right) $ is defined through%
\begin{equation*}
T_{n,X}:=\widetilde{\Phi }\left( \theta ,\mathfrak{P}_{n,X}^{W}\right)
=T\left( \mathfrak{P}_{n,X}^{W}\right)
\end{equation*}%
which allows to reject H0 for large values since $\lim_{n\rightarrow \infty
} $ $T_{n,X}=0$ whenever H0 holds. In the above display we have emphasized
in $\mathfrak{P}_{n,X}^{W}$ the fact that we have used the RV $X_{i}$'s. Let%
\begin{equation*}
L_{n}(t):=P^{W}(\left. T_{n,X}>t\right\vert X_{1},..,X_{n}).
\end{equation*}

We use $P^{W}$ to emphasize the fact that the hazard is due to the weights.\
Consider now a set of RV's $Z_{1},..,Z_{n}$ extracted from a sequence such
that $\lim_{n\rightarrow \infty }P_{n,Z}=P_{\theta ^{\prime }}$ a.s ; we
have denoted $P_{n,Z}$ the empirical measure of $\left(
Z_{1},..,Z_{n}\right) ;$ accordingly define $\mathfrak{P}_{n,Z}^{W^{\prime
}},$ the normalized weighted empirical measure of the $Z_{i}$ 's making use
of weights $\left( W_{1}^{\prime },..,W_{n}^{\prime }\right) $ which are iid
copies of $\left( W_{1},..,W_{n}\right) $, drawn independently from $\left(
W_{1},..,W_{n}\right) .$ Define accordingly 
\begin{equation*}
T_{n,Z}:=\widetilde{\Phi }\left( \theta ,\mathfrak{P}_{n,Z}^{W^{\prime
^{\prime }}}\right) =T\left( \mathfrak{P}_{n,Z}^{W^{\prime ^{\prime
}}}\right) .
\end{equation*}%
Define 
\begin{equation*}
L_{n}(T_{n,Z}):=P^{W}(\left. T_{n,W}>T_{n,Z}\right\vert X_{1},..,X_{n})
\end{equation*}%
which is a RV (as a function of $T_{n,Z}$) . It holds

\begin{equation*}
\lim_{n\rightarrow \infty }T_{n,Z}=\widetilde{\Phi }\left( \theta ,\theta
^{\prime }\right) \text{ \ \ \ a.s}
\end{equation*}%
and therefore the Bahadur slope for the test with statistics $T_{n}$ is $%
\Phi \left( \theta ^{\prime },\theta \right) $ as follows from 
\begin{align*}
\lim_{n\rightarrow \infty }\frac{1}{n}\log L_{n}(T_{n,Z})& =-\inf \left\{
\Phi \left( Q,\theta _{T}\right) :\widetilde{\Phi }\left( \theta ,Q\right) >%
\widetilde{\Phi }\left( \theta ,\theta ^{\prime }\right) \right\} \\
& =-\inf \left\{ \Phi \left( Q,\theta _{T}\right) :\Phi \left( Q,\theta
\right) >\Phi \left( \theta ^{\prime },\theta \right) \right\} \\
& =-\Phi \left( \theta ^{\prime },\theta \right)
\end{align*}

\bigskip if $\theta _{T}=\theta .$ Under H0 the rate of decay of the $p-$%
value corresponding to a sampling under H1 is captured through the
divergence $\Phi \left( \theta ^{\prime },\theta \right) .$

Consider now a competitive test statistics $S\left( \mathfrak{P}%
_{n,X}^{W}\right) $ and evaluate its Bahadur slope.\ Similarly as above it
holds, assuming continuity of the functional $S$ on $\mathbb{S}^{K}$ 
\begin{align*}
\lim_{n\rightarrow\infty}\frac{1}{n}\log P^{W}\left( \left. S\left( 
\mathfrak{P}_{n,X}^{W}\right) >S\left( \mathfrak{P}_{n,Z}^{W^{^{\prime}}}%
\right) \right\vert X_{1},..,X_{n}\right) & =-\inf\left\{ \Phi\left(
Q,\theta_{T}\right) :S(Q)>S\left( \theta^{\prime}\right) \right\} \\
& \geq-\Phi\left( \theta^{\prime},\theta_{T}\right)
\end{align*}
as follows from the continuity of $Q\rightarrow\Phi\left( Q,\theta
_{T}\right) .$ Hence the Bahadur slope of the test based on $S\left( 
\mathfrak{P}_{n,X}^{W}\right) $ is larger or equal $\Phi\left(
\theta^{\prime},\theta\right) .$

We have proved that the chances under H0 for the statistics $T_{n,X}$ to
exceed a value obtained under H1 are (asymptotically) less that the
corresponding chances associated with any other statistics based on the same
bootstrapped sample; as such it is most specific on this scale with respect
to any competing ones. Namely

\begin{proposition}
Under the weighted sampling the test statistics $T\left( \mathfrak{P}%
_{n,X}^{W}\right) $ is most efficient among all tests which are empirical
versions of continuous functionals on $\mathbb{S}^{K}.$
\end{proposition}

\bigskip

\bigskip

\section{Appendix}

\subsubsection{A heuristic derivation of the conditional LDP for the
normalized weighted empirical measure}

The following sketch of proof gives the core argument which\ yields to
Proposition \ref{THM9}; a proof adapted to a more abstract setting can be
found in \cite{WiTr2014}, following their Theorem 3.2 and Corollary 3.3, but
we find it useful to present a proof which reduces to simple arguments. We
look at the probability of the event%
\begin{equation}
P_{n}^{W}\in V(R)  \label{2}
\end{equation}%
for a given vector $R$ in $\mathbb{R}^{K}$, where $V(R)$ denotes a
neighborhood of $R$, therefore defined through%
\begin{equation*}
\left( Q\in V(R)\right) \Longleftrightarrow \left( Q(d_{l})\thickapprox
R(d_{l});1\leq l\leq k\right)
\end{equation*}

We denote by $P$ the distribution of the RV $X$ so that $P_{n}$ $\ $%
converges to $P$ a.s.

Evaluating loosely the probability of the event defined in (\ref{2}) yields,
denoting $P_{X_{1}^{n}}$ the conditional distribution given $\left(
X_{1},..,X_{n}\right) $ 
\begin{align*}
P_{X_{1}^{n}}\left( P_{n}^{W}\in V(R)\right) & =P_{X_{1}^{n}}\left( \bigcap
\limits_{l=1}^{K}\left( \frac{1}{n}\sum_{i=1}^{n}W_{i}\delta_{X_{i}}\left(
d_{l}\right) \thickapprox R(d_{l})\right) \right) \\
& =P_{X_{1}^{n}}\left( \bigcap \limits_{l=1}^{K}\left( \frac{1}{n}%
\sum_{i=1}^{n_{l}}W_{i,l}\thickapprox R(d_{l})\right) \right) \\
& =\prod \limits_{l=1}^{K}P_{X_{1}^{n}}\left( \frac{1}{n_{l}}%
\sum_{i=1}^{n_{l}}W_{i,l}\thickapprox \frac{n}{n_{l}}R(d_{l})\right) \\
& =\prod \limits_{l=1}^{K}P_{X_{1}^{n}}\left( \frac{1}{n_{l}}%
\sum_{i=1}^{n_{l}}W_{i,l}\thickapprox \frac{R(d_{l})}{P(d_{l})}\right)
\end{align*}
where we used repeatedly the fact that the r.v's $W$ are i.i.d. In the above
display, from the second line on, the r.v's are independent copies of $W_{1}$
for all $i$ and $l.$ In the above displays $n_{l}$ is the number of $X_{i}$%
's which equal $d_{l}$, and the $W_{i,l}$ are the weights corresponding to
these $X_{i}$'s. Note that we used the convergence of $n_{l}/n$ to $P(d_{l})$
in the last display.

Now for each $l$ in $\left\{ 1,2,..,K\right\} $ by the Cramer LDP for the
empirical mean, it holds 
\begin{equation*}
\frac{1}{n_{l}}\log P\left( \frac{1}{n_{l}}\sum_{i=1}^{n_{l}}W_{i,l}%
\thickapprox\frac{R(d_{l})}{P(d_{l})}\right) \thickapprox-\varphi^{W}\left( 
\frac{R(d_{l})}{P(d_{l})}\right)
\end{equation*}
i.e.%
\begin{equation*}
\frac{1}{n}\log P\left( \frac{1}{n_{l}}\sum_{i=1}^{n_{l}}W_{i,l}\thickapprox%
\frac{R(l)}{P(l)}\right) \thickapprox-\frac{R(d_{l})}{P(d_{l})}%
\varphi^{W}\left( \frac{R(d_{l})}{P(d_{l})}\right)
\end{equation*}
as follows from the classical Cramer LDP, and therefore 
\begin{align*}
& \frac{1}{n}\log P_{X_{1}^{n}}\left( P_{n}^{W}\in V(R)\right) \\
& \thickapprox\frac{1}{n}\log P_{X_{1}^{n}}\left( \bigcap
\limits_{l=1}^{K}\left( \frac{1}{n}\sum_{i=1}^{n_{l}}W_{i,l}\thickapprox
R(d_{l})\right) \right) \\
& \rightarrow-\sum_{l=1}^{K}\varphi^{W}\left( \frac{R(d_{l})}{P(d_{l})}%
\right) P(d_{l})=-\phi^{W}\left( R,P\right)
\end{align*}

as $n\rightarrow\infty$ .

A precise derivation of Proposition \ref{THM9} involves two arguments:
firstly for a set $\Omega$ $\subset\mathbb{R}^{K}$ a covering procedure by
small balls allowing to use the above derivation locally, and the regularity
assumption (\ref{Omega regul}) which allows to obtain proper limits in the
standard LDP statement.

The argument leading from Proposition \ref{THM9} to Theorem \ref{THM12} can
\ be summarized now.

For some subset $\Omega$ in $\mathbb{S}^{K}$ with non void interior\bigskip
\ it holds%
\begin{equation*}
\left( \mathfrak{P}_{n}^{W}\in\Omega\right) =\bigcup \limits_{m\neq0}\left(
\left( P_{n}^{W}\in m\Omega\right) \cap\left( \sum_{i=1}^{n}W_{i}=m\right)
\right)
\end{equation*}
and $\left( P_{n}^{W}\in m\Omega\right) \subset\left(
\sum_{i=1}^{n}W_{i}=m\right) $ for all $m\neq0.$ Therefore 
\begin{equation*}
P_{X_{1}^{n}}\left( \mathfrak{P}_{n}^{W}\in\Omega\right)
=P_{X_{1}^{n}}\left( \bigcup \limits_{m\neq0}\left( P_{n}^{W}\in
m\Omega\right) \right) .
\end{equation*}
Making use of Theorem \ref{THM9}%
\begin{equation*}
\lim_{n\rightarrow\infty}\frac{1}{n}\log P_{X_{1}^{n}}\left( \mathfrak{P}%
_{n}^{W}\in\Omega\right) =-\phi^{W}\left( \bigcup
\limits_{m\neq0}m\Omega,P\right) .
\end{equation*}
Now 
\begin{equation*}
\phi^{W}\left( \bigcup \limits_{m\neq0}m\Omega,P\right)
=\inf_{m\neq0}\inf_{Q\in\Omega}\phi^{W}\left( mQ,P\right) .
\end{equation*}
We have sketched the arguments leading to Theorem \ref{THM12}; see \cite%
{BrSt2020}\ for details.

\end{document}